\documentclass[12pt]{amsart}
\usepackage{amssymb}
\theoremstyle{plain}
 \newtheorem{thm}{Theorem}[section]
 \newtheorem{cor}[thm]{Corollary}
 \newtheorem{lem}[thm]{Lemma}

\theoremstyle{definition}
 \newtheorem{ex}{Example}[section]
\theoremstyle{remark}
 \newtheorem{rem}{Remark}[section]
\begin{document}
\title[a note on Jensen inequality for 
self-adjoint operators]
{a note on Jensen inequality for self-adjoint operators}
\author[
Tomohiro Hayashi]{{Tomohiro Hayashi} }
\address[Tomohiro Hayashi]
{Nagoya Institute of Technology, 
Gokiso-cho, Showa-ku, Nagoya, Aichi, 466-8555, Japan}
\email[Tomohiro Hayashi]{hayashi.tomohiro@nitech.ac.jp}

\baselineskip=17pt
%\subjclass{47A63, 47A64}
\maketitle
\begin{abstract}
In this paper we consider the order-like relation 
for self-adjoint operators on some Hilbert space. This 
relation is defined by using Jensen inequality. 
We will show that under some assumptions 
this relation is antisymmetric.

\end{abstract}

\section{Introduction}
Let $f(t)$ be a continuous, increasing concave function 
on the real line ${\Bbb{R}}$ and let $A$ be a 
bounded 
self-adjoint operator on some Hilbert space 
${\frak H}$ with an inner product 
$\langle\cdot,\cdot\rangle$. 
Then for each unit vector 
$\xi\in {\frak H}$, we have 
so-called Jensen inequality:
$$\langle 
f(A)\xi,\xi
\rangle
\leq 
f(\langle 
A\xi,\xi
\rangle).
$$
For two self-adjoint operators $X$ and $Y$, if they satisfy 
$
f(X)\leq f(Y)
$, then by using Jensen inequality we have 
$$
\langle 
f(X)\xi,\xi
\rangle
\leq 
\langle 
f(Y)\xi,\xi
\rangle
\leq
f(\langle 
Y\xi,\xi
\rangle). 
$$
Therefore if 
$
\langle 
f(X)\xi,\xi
\rangle
\leq
f(\langle 
Y\xi,\xi
\rangle)
$ for any unit vector 
$\xi\in {\frak H}$, 
we may consider that 
$X$ is dominated by $Y$ 
in some sense. 
Keeping this in our minds, we shall consider the following 
problem: If we have $
\langle 
f(X)\xi,\xi
\rangle
\leq
f(\langle 
Y\xi,\xi
\rangle)
$ and 
$
\langle 
f(Y)\xi,\xi
\rangle
\leq
f(\langle 
X\xi,\xi
\rangle)
$ 
for any unit vector $\xi\in {\frak H}$, 
can we conclude $X=Y$? 
(This problem was suggested by Professor Bourin~\cite{B}.) 

The main results of this paper consist of two 
theorems. 
In section 2 we will solve 
the above problem affirmatively when the 
Hilbert space ${\frak H}$ is finite dimensional. 
Unfortunately we cannot show this 
in the infinite dimensional case. 
But in section 3 we will solve 
a modified problem in full 
generality. 

Here we remark that 
in the paper \cite{A}, 
T.~Ando considered similar problem and 
showed the following theorem: 
``Let $f(t)$ be an operator monotone function. 
If two positive invertible operators $X$ and $Y$ 
satisfy 
$
\langle 
f(X)\xi,\xi
\rangle
\leq
f(\langle 
Y\xi,\xi
\rangle)
$ and 
$
f(\langle 
Y^{-1}\xi,\xi
\rangle^{-1})
\leq 
\langle 
f(X)^{-1}\xi,\xi
\rangle^{-1}
$ 
for any unit vector $\xi\in {\frak H}$, 
then we have $f(X)=f(Y)$.'' 

The author wishes to express his hearty gratitude to Professor 
Jean-Christophe Bourin. 
The author is also grateful to Professor 
Yoshihiro Nakamura for discussion. 
The author would like to thank Professor 
Tsuyoshi Ando for valuable comments. 

\ \\

Throughout this paper we assume that 
the readers are familiar with 
basic notations and results on operator 
theory. We refer the readers to 
Conway's book~\cite{C}. 

We denote 
by ${\frak H}$ a 
(finite or infinite dimensional) 
complex Hilbert space 
and by $B({\frak H})$ 
all bounded linear operators on it. 
For each operator $A\in B({\frak H})$, 
its operator norm is denoted by 
$||A||$. 
For two vectors $\xi,\eta\in {\frak H}$, 
their inner product and norm are denoted by 
$\langle \xi,\eta\rangle$ and $||\xi||$ 
respectively. 
For an interval $[a,b)$, we denote 
its defining function 
by 
$\chi_{[a,b)}(t)$. 
\section{Finite dimensional case}

\begin{thm}
For two hermitian matrices 
$X,Y\in M_{n}({\Bbb C})$ 
and a continuous 
strictly increasing (or decreasing) 
convex function $f(t)$ on 
some interval $I$ containing the 
numerical ranges of $X$ and $Y$, if 
they satisfy 
$$\langle 
f(X)\xi,\xi
\rangle
\geq 
f(\langle 
Y\xi,\xi
\rangle)
$$ 
and 
$$\langle 
f(Y)\xi,\xi
\rangle
\geq 
f(\langle 
X\xi,\xi
\rangle)
$$  
for any unit vector 
$\xi\in {\Bbb C}^{n}$, 
then we have $X=Y$. 

\end{thm}

\begin{proof}
Replacing $f(t)$ by $f(t)+c$ 
for some positive constant $c$ 
if necessarily, we
may assume that $f\geq 0$ on $I$. 
Then $f(X)$ and $f(Y)$ are positive semidefinite matrices. 
Take minimal projections $P$ and $Q$ such that 
$XP=PX$, $YQ=QY$
$f(X)P=||f(X)||P$ and 
$f(Y)Q=||f(Y)||Q$. 
Then for each unit vector $\xi\in Q{\Bbb C}^{n}$ 
we see that 
$\langle 
f(X)\xi,\xi
\rangle Q
=Qf(X)Q
$ 
and 
$
f(\langle 
Y\xi,\xi
\rangle)Q
=||f(Y)||Q
$. Therefore by assumption we have 
$
Qf(X)Q\geq ||f(Y)||Q
$ and hence 
$
||f(X)||Q\geq Qf(X)Q\geq ||f(Y)||Q
$. 
By the similar way we see that 
$
||f(Y)||P\geq Pf(Y)P\geq ||f(X)||P
$. Hence we get $||f(X)||=||f(Y)||$ 
and 
$Qf(X)Q=||f(X)||Q$. 
Since 
$$
0=Q(||f(X)||-f(X))Q
=Q(||f(X)||-f(X))^{\frac{1}{2}}
(||f(X)||-f(X))^{\frac{1}{2}}Q,
$$ 
we have 
$$Qf(X)=f(X)Q=||f(X)||Q=||f(Y)||Q=f(Y)Q$$ 
and hence 
$QX=XQ=YQ$. 
(Here we use the existence of $f^{-1}(t)$.) 
Since two matrices $X(1-Q)$ and $Y(1-Q)$ satisfy 
same assumptions on $(1-Q){\Bbb C}^{n}$, 
we can repeat this argument. Therefore we get 
$X=Y$.

\end{proof}

\begin{cor}
For two hermitian matrices 
$X,Y\in M_{n}({\Bbb C})$ 
and a continuous 
strictly increasing (or decreasing) 
concave function $f(t)$ on some interval 
$I$ containing the 
numerical ranges of $X$ and $Y$, if 
they satisfy 
$$\langle 
f(X)\xi,\xi
\rangle
\leq 
f(\langle 
Y\xi,\xi
\rangle)
$$ 
and 
$$\langle 
f(Y)\xi,\xi
\rangle
\leq 
f(\langle 
X\xi,\xi
\rangle)
$$  
for any unit vector 
$\xi\in {\Bbb C}^{n}$, 
then we have $X=Y$. 

\end{cor}

\begin{proof}
Apply the previous theorem to the function 
$-f(t)$.

\end{proof}

\begin{rem}
If $f(X)$ and $f(Y)$ are of the forms 
$$
f(X)=\sum_{i=1}^{\infty}\lambda_{i}P_{i}
\ \ \ \ \ \ 
f(Y)=\sum_{j=1}^{\infty}\mu_{j}Q_{j}
$$
where $\{P_{i}\}_{i}$ and $\{Q_{j}\}_{j}$ 
are orthogonal family of projections and 
$\lambda_{1}\geq \lambda_{2}\geq\cdots$ 
and $\mu_{1}\geq \mu_{2}\geq\cdots$, 
then Theorem 2.1 holds by the same proof. 
For example, if both $X$ and $Y$ are compact 
positive and $f(t)$ is strictly increasing, 
then $f(X)$ and $f(Y)$ are of the above forms. 

\end{rem}

\section{Infinite dimensional case}
Let $f(t)$ and $g(t)$ be positive, strictly 
increasing, concave $C^{2}$-functions 
on $(0,\infty)$ and continuous 
on $[0,\infty)$. 
For a positive operator $A$, 
by Jensen inequality 
we have 
$$
\langle 
(g\circ f)(A)\xi,\xi
\rangle
\leq 
g(\langle 
f(A)\xi,\xi
\rangle)
\leq 
(g\circ f)(\langle 
A\xi,\xi
\rangle)
$$
for any unit vector 
$\xi\in {\frak H}$. 
We would like to consider the 
``converse'' of this fact. 

\begin{thm}
Let $f(t)$ and $g(t)$ be positive, strictly 
increasing, concave $C^{2}$-functions 
on $(0,\infty)$ and continuous 
on $[0,\infty)$. 
For two positive operators 
$X,Y\in B({\frak H})$, if 
they satisfy 
$$
\langle 
(g\circ f)(X)\xi,\xi
\rangle
\leq 
g(\langle 
f(Y)\xi,\xi
\rangle)
\leq 
(g\circ f)(\langle 
X\xi,\xi
\rangle)
$$
for any unit vector 
$\xi\in {\frak H}$, 
then we have $X=Y$. 

\end{thm}

For example consider the case $f(t)=g(t)=\sqrt{t}$. 
Then we have; 
\begin{ex}
For two positive operators 
$X,Y\in B({\frak H})$, if 
they satisfy 
$$
\langle 
X^{\frac{1}{4}}\xi,\xi
\rangle
\leq 
\langle 
Y^{\frac{1}{2}}\xi,\xi
\rangle^{\frac{1}{2}}
\leq 
\langle 
X\xi,\xi
\rangle^{\frac{1}{4}}
$$
for any unit vector 
$\xi\in {\frak H}$, 
then we have $X=Y$
\end{ex}

The strategy of the proof is essentially same 
as that of \cite{A}\cite{H}. 

\begin{lem}[Ando~\cite{A}]
Let $h(t)$ 
be a positive, strictly 
increasing, concave $C^{2}$-function 
on $(0,\infty)$ and continuous 
on $[0,\infty)$. 
For positive operators $A$ and $B$, 
the inequality 
$$
\langle 
h(A)\xi,\xi
\rangle
\leq 
h(\langle 
B\xi,\xi
\rangle)
$$
holds 
for any unit vector 
$\xi\in {\frak H}$ if and only if we have 
$$
h(A)
\leq 
h'(\lambda)B-\lambda h'(\lambda)+h(\lambda)
$$ 
for any positive number $\lambda$. 

\end{lem}

\begin{proof}
First we will show the 
``only if'' part. 
Since 
$h(t)$ is concave, we have 
$$
h(t)
\leq 
h'(\lambda)t-\lambda h'(\lambda)+h(\lambda).
$$ 
(The right-hand side is the 
tangent line of $h(t)$ at $t=\lambda$.) 
Letting $t=\langle 
B\xi,\xi
\rangle$, we get 
$$
h(\langle B\xi,\xi\rangle)
\leq 
h'(\lambda)\langle 
B\xi,\xi
\rangle-\lambda h'(\lambda)+h(\lambda)
=\langle
\{h'(\lambda)B-\lambda h'(\lambda)+h(\lambda)\}\xi,\xi
\rangle.
$$ 
Combining this with the inequality $\langle 
h(A)\xi,\xi
\rangle
\leq 
h(\langle 
B\xi,\xi
\rangle)
$, we see that 
$$h(A)
\leq 
h'(\lambda)B-\lambda h'(\lambda)+h(\lambda). 
$$ 

Conversely if 
$$h(A)
\leq 
h'(\lambda)B-\lambda h'(\lambda)+h(\lambda) 
$$ holds 
for any $\lambda> 0$, we see that for any unit vector 
$\xi\in {\frak H}$ 
$$
\langle h(A)\xi,\xi\rangle
\leq \langle 
(h'(\lambda)B-\lambda h'(\lambda)+h(\lambda))
\xi,\xi\rangle
=h'(\lambda)\langle 
B\xi,\xi\rangle-\lambda h'(\lambda)+h(\lambda). 
$$ 
Then it is easy to see that the minimal value of 
the right-hand side with respect to 
$\lambda> 0$ is equal to $h(\langle 
B\xi,\xi
\rangle)$.
\end{proof}

\begin{lem}
Under the assumptions in Theorem 3.1, 
we have 
\begin{align*}
\frac{(g\circ f)(X)+f(\lambda) g'(f(\lambda))
-g(f(\lambda))}
{g'(f(\lambda))}
&\leq 
f(Y)\\
&\leq 
f'(\lambda)X-
\lambda f'(\lambda)
+f(\lambda)
\end{align*}
for any positive number $\lambda$. 
\end{lem}

\begin{proof}
By assumptions we have two inequalities 
$$
\langle 
g(f(X))\xi,\xi
\rangle
\leq 
g(\langle 
f(Y)\xi,\xi
\rangle)
$$
and 
$$
\langle 
f(Y)\xi,\xi
\rangle
\leq 
f(\langle 
X\xi,\xi
\rangle)
$$ 
for any unit vector 
$\xi\in {\frak H}$. 
So by the previous lemma we get 
$$
g(f(X))
\leq 
g'(\mu)f(Y)-\mu g'(\mu)+g(\mu)
$$
and 
$$
f(Y)
\leq 
f'(\lambda)X-\lambda f'(\lambda)+f(\lambda).
$$ 
for any positive numbers $\mu$ 
and $\lambda$. Letting $\mu=f(\lambda)$ 
we get the desired inequality. 
\end{proof}

\begin{lem}
Fix two positive numbers $0<a<b$. 
Then there exists a positive constant $c$ 
(depending on the choice of $a,b$) 
such that 
$$
f'(\lambda)t-
\lambda f'(\lambda)
+f(\lambda)-
\left\{
\frac{(g\circ f)(t)+f(\lambda) g'(f(\lambda))
-g(f(\lambda))}
{g'(f(\lambda))}\right\}
\leq c(t-\lambda)^{2}
$$
for any $a\leq\lambda\leq b$ 
and $a\leq t\leq b$.
\end{lem}

\begin{proof}
Set 
$$
k(t)=k_{\lambda}(t)=
c(t-\lambda)^{2}-
f'(\lambda)t+
\lambda f'(\lambda)
-f(\lambda)+
\left\{
\frac{(g\circ f)(t)+f(\lambda) g'(f(\lambda))
-g(f(\lambda))}
{g'(f(\lambda))}
\right\}. 
$$ We will choose an appropriate 
constant $c$ later. Fix $\lambda$ and 
we consider $k(t)$ as a one variable function. 
Then we see that 
$$
k'(t)=
2c(t-\lambda)-
f'(\lambda)+
\frac{(g'\circ f)(t)f'(t)}
{g'(f(\lambda))}
$$
and 
$$
k''(t)=
2c+
\frac{(g''\circ f)(t)f'(t)^{2}+
(g'\circ f)(t)f''(t)}
{g'(f(\lambda))}.
$$
By assumptions we can take $c$ such that $k''(t)>0$ 
for any $a\leq\lambda\leq b$ 
and $a\leq t\leq b$. Then since 
$k'(\lambda)=0$, we have 
$k'(t)\leq 0$ $(t\leq\lambda)$ 
and 
$k'(t)\geq 0$ $(t\geq\lambda)$. 
Hence we have 
$
k(t)\geq k(\lambda)=0.
$
\end{proof}

Take two positive numbers $0<a<b$ such that 
$||X||<b$ 
and $||Y||<b$. 
We can find a 
positive number $\alpha$ 
(depending on the choice of $a,b$) 
such that 
$$
\frac{(g\circ f)(t)+f(\lambda) g'(f(\lambda))
-g(f(\lambda))}
{g'(f(\lambda))}+\alpha\geq 1
$$
for any $a\leq\lambda\leq b$ 
and $a\leq t\leq b$.

\begin{lem}
There exists a positive constant $c$ such that 
$$
\left\{
\frac{(g\circ f)(t)+f(\lambda) g'(f(\lambda))
-g(f(\lambda))}
{g'(f(\lambda))}+\alpha\right\}^{-1}
-
\{f'(\lambda)t-
\lambda f'(\lambda)
+f(\lambda)+\alpha\}^{-1}
\leq c(t-\lambda)^{2}
$$
for any $a\leq\lambda\leq b$ 
and $a\leq t\leq b$. 
The constant $c$ is same as that of the previous 
lemma. 
\end{lem}

\begin{proof}
Set 
$$
p(t)=f'(\lambda)t-
\lambda f'(\lambda)
+f(\lambda)+\alpha
$$
and
$$
q(t)=
\frac{(g\circ f)(t)+f(\lambda) g'(f(\lambda))
-g(f(\lambda))}
{g'(f(\lambda))}+\alpha. 
$$ 
Fix $\lambda$ and 
we consider $p(t)$, $q(t)$ as one variable functions. 
Then $p(t)\geq q(t)\geq 1$ and by the previous lemma 
we have 
$
p(t)-q(t)\leq c(t-\lambda)^{2}.
$ 
So we get 
$$
q(t)^{-1}-p(t)^{-1}=q(t)^{-1}p(t)^{-1}
(p(t)-q(t))\leq c(t-\lambda)^{2}.
$$ 
\end{proof}

\begin{proof}[Proof of Theorem 3.1.]
Take a spectral projection $P$ of $X$. 
By lemma 3.3 we have 
\begin{align*}
\left\{
\frac{(g\circ f)(X)+f(\lambda) g'(f(\lambda))
-g(f(\lambda))}
{g'(f(\lambda))}+\alpha\right\}P
&\leq 
P(f(Y)+\alpha)P\\
&\leq \{(
f'(\lambda)X-
\lambda f'(\lambda)
+f(\lambda))+\alpha \}P
\end{align*}
for any positive number $\lambda$. 
On the other hand we have 
\begin{align*}
\left\{
\frac{(g\circ f)(X)+f(\lambda) g'(f(\lambda))
-g(f(\lambda))}
{g'(f(\lambda))}+\alpha\right\}P
&\leq 
(f(X)+\alpha)P\\
&\leq \{(
f'(\lambda)X-
\lambda f'(\lambda)
+f(\lambda))+\alpha \}P
\end{align*}
for any positive number $\lambda$. 
Combining these with with lemma 3.4 we get 
$$
||
(f(X)+\alpha)P-P(f(Y)+\alpha)P
||
\leq 
c
||
XP-\lambda P
||^{2}\eqno{(1)}
$$ 
whenever 
$P\leq \chi_{[a,b)}(X)$ 
and $a\leq\lambda\leq b$. 

Similarly since we have two inequalities 
\begin{align*}
\{(
f'(\lambda)X-
\lambda f'(\lambda)
+f(\lambda))+\alpha \}^{-1}P
&\leq 
P(f(Y)+\alpha)^{-1}P\\
&\leq\left\{
\frac{(g\circ f)(X)+f(\lambda) g'(f(\lambda))
-g(f(\lambda))}
{g'(f(\lambda))}+\alpha\right\}^{-1}P
\end{align*}
and 
\begin{align*}
\{(
f'(\lambda)X-
\lambda f'(\lambda)
+f(\lambda))+\alpha \}^{-1}P
&\leq 
(f(X)+\alpha)^{-1}P\\
&\leq\left\{
\frac{(g\circ f)(X)+f(\lambda) g'(f(\lambda))
-g(f(\lambda))}
{g'(f(\lambda))}+\alpha\right\}^{-1}P, 
\end{align*}
by lemma 3.5 we get 
$$
||
(f(X)+\alpha)^{-1}P-P(f(Y)+\alpha)^{-1}P
||
\leq 
c
||
XP-\lambda P
||^{2}
$$ 
whenever 
$P\leq \chi_{[a,b)}(X)$ 
and $a\leq\lambda\leq b$. 
Hence in this case 
\begin{align*}
||
(f&(X)+\alpha)P-(P(f(Y)+\alpha)^{-1}P)^{-1}
||\\
&=
||
(f(X)+\alpha)
\{
P(f(Y)+\alpha)^{-1}P
-(f(X)+\alpha)^{-1}P
\}(P(f(Y)+\alpha)^{-1}P)^{-1}
||\\
&\leq 
||
f(X)+\alpha||\cdot ||
(P(f(Y)+\alpha)^{-1}P)^{-1}||\cdot||
(P(f(Y)+\alpha)^{-1}P
-(f(X)+\alpha)^{-1}P
||\\
&\leq 
(f(b)+\alpha)^{2}c
||
XP-\lambda P
||^{2}.
\end{align*}
Therefore for $P\leq \chi_{[a,b)}(X)$ 
and $a\leq\lambda\leq b$ we have 
$$
||
P(f(Y)+\alpha)P
-
(P(f(Y)+\alpha)^{-1}P)^{-1}
||
\leq (1+(f(b)+\alpha)^{2})c
||
XP-\lambda P
||^{2}.\eqno{(2)}
$$
The rest of the proof is almost same as that 
of \cite{A}\cite{H}. We include this for the 
reader's convenience. 

For each integer $n$, let $P_{i}\ \ (i=1,2,\cdots,n)$ 
be the spectral projections of $X$ corresponding to 
the interval 
$[a+\frac{(i-1)(b-a)}{n},a+\frac{i(b-a)}{n})$. 
Then we have $\sum_{i}P_{i}=\chi_{[a,b)}(X)$ and 
$$
||
XP_{i}-\lambda_{i}P_{i}
||
\leq 
\frac{b-a}{n}
$$
where $\lambda_{i}=a+\frac{(i-1)(b-a)}{n}$. 
Then it follows from (1) that 
$$
||\sum_{i=1}^{n}\{
(f(X)+\alpha)P_{i}-P_{i}(f(Y)+\alpha)P_{i}\}
||\leq 
\frac{c(b-a)^{2}}{n^{2}}.\eqno{(3)}
$$
Similarly it follows from (2) that 
$$
||P_{i}(f(Y)+\alpha)P_{i}
-
(P_{i}(f(Y)+\alpha)^{-1}P_{i})^{-1}||
\leq 
\frac{(1+(f(b)+\alpha)^{2})c(b-a)^{2}}{n^{2}}.
$$
By using the 
following formula, which is so-called 
Schur complement 
$$
(P_{i}(f(Y)+\alpha)^{-1}P_{i})^{-1}
=P_{i}(f(Y)+\alpha)P_{i}
-P_{i}(f(Y)+\alpha)P_{i}^{\perp}
(P_{i}^{\perp}(f(Y)+\alpha)P_{i}^{\perp})^{-1}
P_{i}^{\perp}(f(Y)+\alpha)P_{i}
$$
where $P_{i}^{\perp}=1-P_{i}$, 
we see that 
\begin{align*}
||P_{i}^{\perp}(f(Y)+\alpha)P_{i}||^{2}
&=||(P_{i}^{\perp}(f(Y)+\alpha)P_{i}^{\perp})^{1/2}
(P_{i}^{\perp}(f(Y)+\alpha)P_{i}^{\perp})^{-1/2}
P_{i}^{\perp}(f(Y)+\alpha)P_{i}||^{2}\\
&\leq 
||f(Y)+\alpha||\cdot
||(P_{i}^{\perp}(f(Y)+\alpha)P_{i}^{\perp})^{-1/2}
P_{i}^{\perp}(f(Y)+\alpha)P_{i}||^{2}\\
&=||f(Y)+\alpha||\cdot
||P_{i}(f(Y)+\alpha)P_{i}^{\perp}
(P_{i}^{\perp}(f(Y)+\alpha)P_{i}^{\perp})^{-1}
P_{i}^{\perp}(f(Y)+\alpha)P_{i}||\\
&=
||f(Y)+\alpha||\cdot
||P_{i}(f(Y)+\alpha)P_{i}-
(P_{i}(f(Y)+\alpha)^{-1}P_{i})^{-1}||\\
&\leq 
\frac{(f(b)+\alpha)(1+(f(b)+\alpha)^{2})
c(b-a)^{2}}{n^{2}}. 
\end{align*}
Therefore by the well-known formula 
$||A||^{2}=||AA^{*}||=||A^{*}A||$ 
we see that 
\begin{align*}
||\sum_{i=1}^{n}
P_{i}^{\perp}(f(Y)+\alpha)P_{i}||^{2}
&=||\{\sum_{i=1}^{n}
P_{i}^{\perp}(f(Y)+\alpha)P_{i}\}\{
\sum_{j=1}^{n}P_{j}(f(Y)+\alpha)P_{j}^{\perp}\}||\\
&=
||\sum_{i=1}^{n}
P_{i}^{\perp}(f(Y)+\alpha)P_{i}(f(Y)+\alpha)P_{i}^{\perp}||\\
&\leq 
\sum_{i=1}^{n}
||P_{i}^{\perp}(f(Y)+\alpha)P_{i}(f(Y)+\alpha)P_{i}^{\perp}||\\
&=\sum_{i=1}^{n}
||P_{i}^{\perp}(f(Y)+\alpha)P_{i}||^{2}\\
&\leq 
\sum_{i=1}^{n}\frac{(f(b)+\alpha)(1+(f(b)+\alpha)^{2})
c(b-a)^{2}}{n^{2}}\\
&=\frac{(f(b)+\alpha)(1+(f(b)+\alpha)^{2})
c(b-a)^{2}}{n}.
\end{align*}
Thus we get 
$$
||\sum_{i=1}^{n}
P_{i}^{\perp}(f(Y)+\alpha)P_{i}||
\leq 
\sqrt{\frac{(f(b)+\alpha)(1+(f(b)+\alpha)^{2})
c(b-a)^{2}}{n}}
.\eqno{(4)}
$$
Since 
$$
f(Y)\chi_{[a,b)}(X)=\sum_{i=1}^{n}P_{i}(f(Y)+\alpha)P_{i}
+\sum_{i=1}^{n}P_{i}^{\perp}(f(Y)+\alpha)P_{i},
$$
by using (3) and (4) we see that 
\begin{align*}
||f(X)&\chi_{[a,b)}(X)-f(Y)\chi_{[a,b)}(X)||\\
&\leq 
||\sum_{i=1}^{n}\{
(f(X)+\alpha)P_{i}
-P_{i}(f(Y)+\alpha)P_{i}\}||
+||\sum_{i=1}^{n}
P_{i}^{\perp}(f(Y)+\alpha)P_{i}||\\
&\leq 
\frac{c(b-a)^{2}}{n^{2}}
+\sqrt{\frac{(f(b)+\alpha)(1+(f(b)+\alpha)^{2})
c(b-a)^{2}}{n}}
.
\end{align*}
By tending $n\rightarrow\infty$ we get 
$f(X)\chi_{[a,b)}(X)=f(Y)\chi_{[a,b)}(X)$. 
Since $a$ is arbitrary we have 
$f(X)\chi_{(0,b)}(X)=f(Y)\chi_{(0,b)}(X)$. 
Therefore in order to show $f(X)=f(Y)$, 
now it is enough to show that 
$\chi_{\{0\}}(X)=\chi_{\{0\}}(Y)$. 

For any unit vector $\xi\in {\frak H}$ such that 
$X\xi=0$, we see that 
$$f(0)+
\langle 
(f(Y)-f(0))\xi,\xi
\rangle=
\langle 
f(Y)\xi,\xi
\rangle
\leq 
f(\langle 
X\xi,\xi
\rangle)=f(0).
$$
Therefore $f(Y)\xi=f(0)\xi$ and hence 
$Y\xi=0$. 
Conversely for any unit vector 
$\xi\in {\frak H}$ such that 
$Y\xi=0$, we see that 
$$(g\circ f)(0)+
\langle 
((g\circ f)(X)-(g\circ f)(0))\xi,\xi
\rangle=
\langle 
(g\circ f)(X)\xi,\xi
\rangle
\leq 
g(\langle 
f(Y)\xi,\xi
\rangle)=(g\circ f)(0).
$$
Therefore 
$(g\circ f)(X)\xi=(g\circ f)(0)\xi$ and hence 
$X\xi=0$. 
\end{proof}

\begin{rem}
\begin{enumerate}
\item
In lemma 3.4, the assumption $a>0$ is 
crucial. For example 
if we consider the case $a=0$ and $f(t)=g(t)=\sqrt{t}$, 
then lemma 3.4 is wrong. Indeed in this case 
\begin{align*}
f'(\lambda)&t-
\lambda f'(\lambda)
+f(\lambda)-
\left\{
\frac{(g\circ f)(t)+f(\lambda) g'(f(\lambda))
-g(f(\lambda))}
{g'(f(\lambda))}\right\}\\
&=
\frac{t}{2\sqrt{\lambda}}
+\frac{3\sqrt{\lambda}}{2}-
2\lambda^{\frac{1}{4}}
t^{\frac{1}{4}}. 
\end{align*}
It is easy to see that 
$$
\frac{1}{(t-\lambda)^{2}}
\left\{
\frac{t}{2\sqrt{\lambda}}
+\frac{3\sqrt{\lambda}}{2}-
2\lambda^{\frac{1}{4}}
t^{\frac{1}{4}}
\right\}
$$ 
is unbounded for 
$0< \lambda\leq b$ and $0< t\leq b$. 
(Fix $t>0$ and consider the case 
$\lambda\rightarrow +0$. Then this function 
tends to $\infty$.)
\item 
The argument in this section cannot be 
applied directly to the problem in the previous section. 
For simplicity, we would like consider the case 
$f(t)=\sqrt{t}$. 
Let $X$ and $Y$ be positive operators 
on ${\frak H}$. 
Suppose that 
they satisfy 
$$\langle 
\sqrt{X}\xi,\xi
\rangle
\leq 
\sqrt{\langle 
Y\xi,\xi
\rangle)}
$$ 
and 
$$\langle 
\sqrt{Y}\xi,\xi
\rangle
\leq 
\sqrt{\langle 
X\xi,\xi
\rangle}
$$  
for any unit vector 
$\xi\in {\frak H}$. Then by lemma 3.2 
we have 
$$
\sqrt{X}\leq \frac{1}{2\sqrt{\lambda}}Y
+\frac{\sqrt{\lambda}}{2}
$$
and 
$$
\sqrt{Y}\leq \frac{1}{2\sqrt{\lambda}}X
+\frac{\sqrt{\lambda}}{2}
$$
for any $\lambda> 0$. 
By the first inequality we have 
$$
2\sqrt{\lambda X}-\lambda
\leq Y.
$$ Since the left-hand side 
in this inequality 
is not positive, 
we cannot take a square root. This is the main trouble. 
By this reason we cannot show the statement 
like lemma3.3.
\end{enumerate}
\end{rem}


\begin{thebibliography}{99}

\bibitem{A} 
T.~Ando, 
{\it Functional calculus with operator-monotone 
functions}, 
Math. Inequal. Appl. (to appear)

\bibitem{B} 
J-C.~Bourin, private 
communication,

\bibitem{C} 
J.~B.~Conway, {\it A course in operator theory.} 
Graduate Studies in Mathematics, 21. 
American Mathematical Society, Providence, RI, 2000.

\bibitem{H}
T.~Hayashi, 
{\it Non-commutative A-G mean inequality.} 
Proc. Amer. Math. Soc. (to appear)

\end{thebibliography}
\end{document}